\documentclass[11pt,a4paper]{article}
 \usepackage{amsmath,amssymb,latexsym,hyperref,xcolor,amsthm,authblk}
 \def\le{\leqslant}
 \def\ge{\geqslant}
 \def\Z{\mathbb Z}
 \def\N{\mathbb N}
  \def\I{\mathbb I}

 \def\epsilon{\varepsilon}
 \def\phi{\varphi}
  \newtheorem{defi}{Definition}
  \newcommand{\bd}{\begin{defi}}      
  \newcommand{\ed}{\end{defi}}

  \newtheorem{lemm}[defi]{Lemma}  
  \newcommand{\bl}{\begin{lemm}}
  \newcommand{\el}{\end{lemm}} 

  \newtheorem{theo}[defi]{Theorem}
  \newcommand{\bt}{\begin{theo}}
  \newcommand{\et}{\end{theo}}

  \newtheorem{cor}[defi]{Corollary}
  \newcommand{\bc}{\begin{cor}}
  \newcommand{\ec}{\end{cor}}

  \newtheorem{pro}[defi]{Proposition}
  \newcommand{\bp}{\begin{pro}}
  \newcommand{\ep}{\end{pro}}
  \theoremstyle{definition}
    \newtheorem{rmk}[defi]{Remark}
  \newcommand{\brmk}{\begin{rmk}}
  \newcommand{\ermk}{\end{rmk}}
  
      \newtheorem{ex}[defi]{Example}
  \newcommand{\bex}{\begin{ex}}
  \newcommand{\eex}{\end{ex}}

  \makeatletter
  \def\proof{\@ifnextchar[\opargproof{\opargproof[\bf Proof \hfil\\ ]}}
  \def\opargproof[#1]{\par\noindent {\bf #1 }}
  \def\endproof{{\unskip\nobreak\hfil\penalty50\hskip8mm\hbox{}
  \nobreak\hfil
  \(\Box\)\parfillskip=0mm \par\vspace{3mm}}}
  \makeatother

\definecolor{red}{rgb}{0, 0, 0}
\definecolor{blue}{rgb}{0, 0, 0}

\begin{document}
\title{\LARGE\bf The Ring of Polyfunctions over $\Z/n\Z$}
\author[1]{Ernst Specker}
\author[2]{Norbert Hungerb\"uhler}
\author[3]{Micha Wasem}
\affil[1]{\em Dedicated to the memory of the first author\bigskip}
\affil[1]{Department of Mathematics, ETH Z\"urich, R\"amistrasse 101, 8092 Z\"urich, Switzerland\smallskip}
\affil[2]{Department of Mathematics, ETH Z\"urich, R\"amistrasse 101, 8092 Z\"urich, Switzerland\\{\tt norbert.hungerbuehler@math.ethz.ch}\smallskip}
\affil[3]{HTA Freiburg, HES-SO University of Applied Sciences and Arts Western Switzerland, P\'erolles 80, 1700 Freiburg, Switzerland\\{\tt micha.wasem@hefr.ch}\smallskip}
\affil[3]{UniDistance Suisse, Schinerstrasse 18, 3900 Brig, Switzerland\\{\tt micha.wasem@fernuni.ch}}

 
\maketitle
 \begin{abstract}
We study the ring of \emph{polyfunctions} over \textcolor{red}{$\Z/n\Z$. The ring of polyfunctions over }a commutative ring $R$ with unit element \textcolor{red}{is} the ring of functions $f:R\to R$ which admit a polynomial representative $p\in R[x]$ in the sense that $f(x)= p(x)$ for all $x\in R$. This allows to define a ring invariant $s$ which
associates to a commutative ring $R$ with unit element a value in $\N\cup\{\infty\}$.  The function $s$ generalizes the number theoretic Smarandache function.
For the ring $R=\Z/n\Z$ we provide a unique representation of polynomials which vanish as a function.
This yields a new formula for the number $\Psi(n)$ of polyfunctions over $\Z/n\Z$.
We also investigate algebraic properties of the ring of polyfunctions over $\Z/n\Z$.
In particular, we identify the additive subgroup of the ring
and the ring structure itself. Moreover we derive formulas for the size of the ring of polyfunctions in several variables
over $\Z/n\Z$, and we compute the number of polyfunctions which are units of the ring.
 \end{abstract}


\section{Introduction}
In a finite field $F$, every function $f: F\to F$ can
be represented by a polynomial, i.e., there exists a polynomial
$p\in F[x]$ such that $f(x)=p(x)$ for all $x\in F$. Such a polynomial
is, e.g., given by the Lagrange interpolation polynomial for $f$.
Among the commutative rings with unit element, the finite 
fields are actually characterized by this representation property
(see~\cite{redei2}): 
\bt[R\'edei, Szele]\label{iff}
If $R$ is a commutative ring with unit element then $R$ is a finite field if and only if every function $f:R\to R$ can be represented by a polynomial in $R[x]$.

\et
If a commutative ring $R$ with unit element
is {\em not\/} a field, it is natural to ask what can be said about the functions
from $R$ to $R$ which {\em can\/} be represented by a polynomial
in $R[x]$. These functions are called polynomial functions or \emph{polyfunctions\/} for short.
The set of polyfunctions
$$\{f:R\to R\mid \exists p\in R[x]\,\, \forall x\in R :  p(x)=f(x)\},$$
equipped with pointwise addition and multiplication,
is a subring of $R^{R}$. This ring of polyfunctions over $R$ will
be denoted by $G(R)$. Of particular interest are the polynomials which correspond to the zero 
element in $G(R)$, they will be called \emph{null-polynomials} (see, e.g.,~\cite{singmaster}).
It is the objective of
this article to investigate  the algebraic structure and combinatorial properties
of the ring of polyfunctions $G(\Z/n\Z)$. 

More generally, one can study the ring of multivariate polyfunctions in $d\in\N$ variables -- this ring is defined as the set
$$\{f:R^d\to R\mid 
\exists p\in R[x_1,x_2,\ldots,x_d]\,\, \forall x=(x_1,\ldots,x_d)\in R^d : p(x)=f(x)\},$$ 
equipped with pointwise addition and multiplication.
We denote this ring by $G_d(R)$ and write $G(R)=G_1(R)$, in accordance with the notation introduced above.

Polyfunctions in one variable over $\Z/n\Z$ were already discussed by Kempner~\cite{kempner2,kempner3}, who gave a formula for the number $\Psi(n)$ of polyfunctions over $\Z/n\Z$, which was subsequently simplified by Keller and Olson in~\cite{keller} (see also the work of Carlitz~\cite{carlitz} in the case where $n$ is a power of a prime). 
Regarding polyfunctions in $d$ variables we refer to Mullen~\cite{mullen} and more recently to~\cite{noebi}: In~\cite[Theorem 2, p.\,5]{noebi}, a characterization theorem is proved which allows to tell whether a given function $f:(\Z/n\Z)^d\to\Z/n\Z$ is a polyfunction or not. 
Furthermore, a formula for the number of polyfunctions $\Psi_d(n)$ in $d$ variables over $\Z/n\Z$ is obtained.
In the present work, we provide an alternative formula for $\Psi(n)$ and a new proof of the formula for $\Psi_d(n)$ given in~\cite{noebi}.

Polyfunctions from $\Z/n\Z$ to $\Z/m\Z$ have been
discussed by Chen~\cite{chen}, \cite{chen2} and~Bhargava~\cite{bhargava}.
The focus there is to find conditions on the pair $(m,n)$
such that all functions (or certain subclasses)
from $\Z/n\Z$ to $\Z/m\Z$ are polyfunctions. These results have been generalized to
polynomial functions in the residue class rings of Dedekind domains by Li and Sha in~\cite{li}.
Dueball in~\cite{dueball} considered polynomials $\bmod$ $p^n$ with  integer coefficients.
He showed that the values of such a polynomial $f(x)$ are already determined when $x$ runs through a certain subset of residues. 
He also provided a formula to generate polynomials which vanish $\bmod$ $p^n$ for all integral values of $x$.

To each commutative ring $R$ with unit element, we can associate a number $s(R)\in\N\cup\{\infty\} $ which is defined to be the minimal degree $m$ such that the function $x\mapsto x^m$ can be represented by a polynomial in $R[x]$ of degree strictly smaller than $m$, i.e.\
\begin{equation}\label{smarandachering}
s(R) := \min\{m\in\N \mid \exists p\in R[x], \operatorname{deg}(p)<m, \forall x\in R : p(x) = x^m\}
\end{equation}
if such an $m$ exists, and $s(R)=\infty$ otherwise.

\textcolor{red}{If $s(R)$ is finite, the monomial $x^{s(R)}$ can be represented by a polynomial $p$ of degree less than $s(R)$. Therefore, the normed
polynomial $q(x) = x^{s(R)}-p(x)$ represents the zero-function. Vice versa, if $r(x)$ is a
normed null-polynomial of minimal degree $m$, then $m=s(R)$.  
Hence, $s(R)$ can be interpreted as the minimal degree of a normed null-polynomial over $R$.}

An alternative and, for reasons that will become clear later,
preferable way to view \textcolor{black}{the function defined by~\eqref{smarandachering}} is as follows: The building blocks
of polynomials are the monomials $x^0,x^1,x^2,\ldots$. We say, a
monomial $x^m$ is {\em reducible}, if the function $x\mapsto x^m$ can be represented by a 
polynomial in $R[x]$ of degree strictly smaller than $m$. 
Then, $s(R)$ is the number of non-reducible monomials.

The function $s$ is a ring invariant which generalizes the classical number theoretic Smarandache function $s:\N\to \N$,
\begin{equation}\label{classical_smarandache}
n\mapsto s(n) :=\min\{k\in\N: n\mid k!\},
\end{equation}
which is named after the \textcolor{red}{Romanian} mathematician
Florentin Smarandache, but which has been originally introduced by Lucas in~\cite{lucas} (for prime powers) and Kempner in~\cite{kempner} (for general $n$). The function $s$ defined in \eqref{smarandachering} will be called {\em Smarandache function} because $n\mapsto s(\Z/n\Z)$
coincides with the usual Smarandache function $n\mapsto s(n)$ (see Theorem~\ref{smarandache} below). In the context of general commutative rings with unit element, this function will be studied in a forthcoming paper~\cite{artikel2}. We also refer to~\cite{redei}, where polyfunctions over general rings are discussed. 

%

The article is organized as follows: 
Section~\ref{z} establishes a unique representation theorem for null-polynomials (Theorem~\ref{deco}).
This provides a new formula for the number $\Psi(n)$ of polyfunctions over $\Z/n\Z$ (Corollary~\ref{co1} and Proposition~\ref{ppp}).
In Section~\ref{algebra}, we investigate algebraic properties of the ring of polyfunctions over $\Z/n\Z$.
In particular, we identify the additive subgroup of the ring (Theorem~\ref{thth})
and the ring structure itself (Theorem~\ref{ring}). We also investigate the multiplicative subgroup $U_n$ of units
in the ring (Propositions~\ref{units} and~\ref{units-general}).
 \textcolor{black}{Section~\ref{xy} comprises
a description of the ring of polyfunctions in several variables
over $\Z/n\Z$. In particular, we give a new formula for the size of this ring (Proposition~\ref{alternative}).} 

\subsection{Notational conventions}
\def\x{{\text{\boldmath$x$}}}\def\k{{\text{\boldmath$k$}}}\def\l{{\text{\boldmath$l$}}}
\textcolor{black}{Unless stated otherwise, $n$ will denote a natural number $\ge 2$ and $\Z_{n}=
\Z/n\Z$ is the ring of integers modulo $n$.} We adopt the notation
$(a,b)$ for the greatest common divisor of the integer numbers $a$ 
and $b$, and we write $a\mid b$ if $b$ is an integer multiple of $a$. \textcolor{red}{Furthermore, for $f,g\in\Z_n[x]$ we will write $f\equiv g \mod n$ to mean the equality of polynomials and we will write $f(x) \equiv g(x) \mod n$ if the functions defined by $f$ and $g$ agree.}

\section{Combinatorial aspects of Polyfunctions over $\Z_{n}$}\label{z}
\subsection{The Smarandache Function}
In this section, we want to determine the minimal degree of
a normed null-polynomial in $\Z_n[x]$. We call a polynomial
{\em normed}, if its leading coefficient is 1.
The answer is given in the following 
theorem:
\bt\label{smarandache}
$s(\Z_n)$ equals the Smarandache function $s(n)$ defined in~\eqref{classical_smarandache}. 
\et
\begin{rmk}
\textcolor{black}{
According to our conventions, $n\ge 2$ as the case $n=1$ should formally be excluded since $\Z_1$ is not a ring with unit element. However, if $n=1$ we can still make sense of $s(\Z_1)$ if we view $\Z_1$ as $\{0\}$ and it holds that $s(\Z_1)=0$ but $s(1) = 1$. Kempner originally defined $s(1)=1$ in~\cite{kempner} but changed it to $s(1) = 0$ later on in~\cite{kempner2,kempner3}. By defining
$$
s(n):=\min\{k\in \N_0:n\mid k!\},
$$
this ambiguity can be avoided (see also~\cite[p.\,7]{noebi}) and the theorem might be stated for every $1\leqslant n\in\N$. Another proof of Theorem~\ref{smarandache} also appears in~\cite[Theorem 7, p.\,126]{noebiundlorenz}.}
\end{rmk}
In order to prove Theorem~\ref{smarandache} for $n\ge 2$, we first show that $s(\Z_n)\le s(n)$. This is established by giving
a normed null-polynomial of degree $s(n)$.  In fact, we have
\begin{align*}
p(x)&:=\prod_{i=1}^{s(n)}(x+i)=\binom{x+s(n)}{s(n)} s(n) ! \,\equiv\, 0 \mod n
\end{align*}
for all $x\in \Z_n$.

The second step consist in proving the reverse inequality $s(\Z_n)\ge s(n)$.
This follows easily from the combinatorial identity
which connects the binomial and the Stirling numbers of the second kind 
(see, e.g.,~\cite[3.39, p.\,97]{aigner} \textcolor{red}{or~\cite[Lemma 3]{noebiundlorenz}}):
For all $r,j\in \N_0$ there holds
$$\sum_{i=0}^r (-1)^{r-i}\binom ri i^j = r ! \left\{ j\atop r\right\}$$
(with the convention $0^0:=1$).
In particular, it follows that
\begin{equation}\label{special}
\sum_{i=0}^r(-1)^{i+r}\binom ri i^k =\delta_{kr} r ! 
\end{equation}
for $k\in\{0,1,\ldots,r\}$. Now, we consider a null-polynomial $p$
over $\Z_n$, i.e., we assume
$$p(i)=\sum_{k=0}^r a_k i^k\equiv 0\mod n$$
for all $i\in\Z_n$. Then, it follows from~(\ref{special}) that modulo $n$
\begin{eqnarray*}
0&\equiv& \sum_{i=0}^r\sum_{k=0}^r (-1)^{i+r}\binom ri a_k i^k\\
 &=     & \sum_{k=0}^ra_k \sum_{i=0}^r (-1)^{i+r}\binom ri i^k\\
 &=     & \sum_{k=0}^ra_k\delta_{kr} r ! =a_r r ! 
\end{eqnarray*}
This establishes the desired inequality $s(\Z_n)\ge s(n)$ 
and the proof of Theorem~\ref{smarandache} is complete.\endproof

In order to gain more insight in the ideal of null-polynomials in $\Z_n[x]$,
we we need a stronger version of Theorem~\ref{smarandache}.
First we consider the following simple lemma:
\bl\label{matrix}
Let $A$ and $C$ denote matrices with integer coefficients, 
$y$ a vector with integer components and $\I$ the identity matrix.
If $A^tC\equiv m\I\mod n$, then $Ay\equiv 0\mod n$
implies $my\equiv 0\mod n$.
\el
\proof Modulo $n$ we have
\begin{align*}
0\equiv C^tAy=(y^tA^tC)^t\equiv (y^tm\I)^t=my.\tag*{$\Box$}
\end{align*}

Lemma~\ref{matrix} allows to prove the following stronger form
of Theorem~\ref{smarandache}. \textcolor{red}{This will be the technical key
to the understanding of the null-polynomials in Section~\ref{dec},
the structure of the additive group of the polyfunctions in Section~\ref{addi},  and of their ring structure
in Section~\ref{riri}.}
\bt\label{voll}
If $p(x)=a_0+a_1 x+\ldots+a_r x^r$ vanishes in $\Z_n$
on the set $x\in\{\alpha,\alpha+1,\ldots,\alpha+r\}$
(in particular, if $p$ is a null-polynomial over $\Z_n$),
then $a_k r ! \equiv 0\bmod n$ holds for all $k\in\{0,1,\ldots,r\}$.
\et
\proof
For $\alpha\in\{0,1,\ldots,n-1\}$ and
$j\in\{\alpha,\alpha+1,\ldots,\alpha+r\}$, we consider the polynomials
$$g_{j,\alpha}(x):=\prod_{\substack{k=\alpha\\ k\neq j}}^{\alpha+r}(x-k) \textcolor{red}{\ =\sum_{k=0}^rg_{j\alpha k}x^k.}$$
Obviously, we have $g_{j,\alpha}(i)=0$ whenever 
$i\in\{\alpha,\alpha+1,\ldots,\alpha+r\}$
is different from $j$, and $g_{j,\alpha}(j)=(j-\alpha) ! 
(-1)^{\alpha+r-j}(\alpha+r-j) ! $.
Hence, we obtain for $i,j\in\{\alpha,\alpha+1,\ldots,\alpha+r\}$ 
$$(-1)^{\alpha+r-j}\binom r{j-\alpha} g_{j,\alpha}(i)=\delta_{ij} r ! $$
\textcolor{red}{This identity can be read as $AD=r! \I$ for the matrix $(A)_{ik}=i^k$, 
$i\in\{\alpha,\alpha+1,\ldots,\alpha+r\}$, $k\in\{0,1,\ldots,r\}$,
and the matrix 
$$
(D)_{kj} = (-1)^{\alpha+r-j}\binom{r}{j-\alpha}g_{j\alpha k},
$$
$k\in\{0,1,\ldots r\}$, $j\in\{\alpha,\alpha+1,\ldots \alpha+r\}$.
Finally, from $AD=r ! \I$ it follows $A^tC= r! \I$ for $C=D^t$.}
Thus, the hypotheses of Lemma~\ref{matrix} are fulfilled
with $m=r ! $.

From the hypothesis of Theorem~\ref{voll} it follows moreover,
that $Ay\equiv 0\mod n$ for the vector $y=(a_0,a_1,\ldots,a_r)^t$
and hence, the conclusion of Lemma~\ref{matrix} gives the
desired result.
\endproof

\subsection{Decomposition of Null-Polynomials}\label{dec}
In this section we analyse the null-polynomials in $\Z_n[x]$,
i.e.\ the polynomials which vanish as a function from 
$\Z_n$ to $\Z_n$. In particular we will determine the number
of null-polynomials which then allows to compute the number of
polyfunctions over $\Z_n$.

We introduce the following notation for $2\leqslant n\in\N$:
$q(n)$ denotes the smallest prime divisor of $n$,
$t(n): = \operatorname{card}\left\{s((n,\alpha!))\mid s((n,\alpha!))\geqslant q(n),\alpha\in\N\right\}
$
and
$$\left\{s((n,\alpha!))\mid s((n,\alpha!))\geqslant q(n),\alpha\in\N\right\}=:\{\beta_1,\beta_2,\ldots, \beta_{t(n)}\},
$$
where the numbers $\beta_k$ are numbered in descending order, i.e.\
\begin{equation}\label{bet}
s(n) = \beta_1 > \beta_2 > \ldots >\beta_{t(n)}=q(n).
\end{equation}
Here, $s$ continues to denote the number-theoretic Smarandache function.
\textcolor{red}{We have
$$\beta_{l+1}=s((n,(\beta_{l}-1)!))$$
for $l=1,\ldots,t(n)-1$: To see this, let $\alpha\in\{q(n),q(n)+1,\ldots,s(n)\}$ be such that $\beta_l = s((n,\alpha!))$.
If $k=(n,\alpha!)$, then $s(k)$ is the smallest number such that $k\mid s(k)!$. If $\alpha>s(k)$ we might replace $\alpha$ by $s(k)$ and obtain $(n,\alpha!)=(n,s(k)!)=(n,\beta_l!)$. Therefore $\beta_l = s((n,\beta_l!))$ and
$
\beta_{l+1} =s((n,(\beta_l-1)!))<\beta_l,
$ as claimed.
}

Furthermore, we define
\begin{equation}\label{alp}
\alpha_k:=\frac{n}{(n,\beta_k ! )}
\end{equation}
and consider the 
{\em basic null-polynomials\/} in $\Z_n[x]:$
\textcolor{red}{
\begin{equation}\label{basicnullpolynomial} b_k(x) := \alpha_k \prod_{i=1}^{\beta_k}(x+i)\end{equation}}

Why the null-polynomials are important becomes clear in Theorem~\ref{deco}
below. But first we consider an example and give some
computational remarks.

\bex \label{counting-zero-polynomials-1}The smallest prime divisor of $n=90$ is $q(90)=2$,
and $s(90)=6$.
In order to compute the degrees $\beta_k$ according to~(\ref{bet}), 
notice that
we only need to consider values  $\alpha\in \{q(n),q(n)+1,
\ldots,s(n)\}$. For these values, we have
\begin{center}
\begin{tabular}{|c|c|c|}\hline
$\alpha$ & $(90,\alpha ! )$ & $s((90,\alpha ! ))$\\\hline\hline
2      & 2              & 2\\
3      & 6              & 3\\
4      & 6              & 3\\
5      & 30             & 5 \\
6      & 90             & 6\\\hline
\end{tabular}
\end{center}
From this table we read off \textcolor{black}{$t(90)=4$} and
$$\beta_1=6,\,\,\beta_2=5,\,\,\beta_3=3,\,\,\beta_4=2.$$
The coefficients $\alpha_k$ are now computed by~(\ref{alp}):
$$\alpha_1=1,\,\,\alpha_2=3,\,\,\alpha_3=15,\,\,\alpha_4=45.$$
The basic null-polynomials for $n=90$ are therefore
\begin{align}
 b_1(x)&=    (1 + x)(2 + x)(3 + x)(4 + x) (5 + x)(6 + x)\notag\\
 b_2(x)&=     3 (1 + x) (2 + x) (3 + x) (4 + x)(5 + x)\notag\\
 b_3(x)&=     15 (1 + x) (2 + x) (3 + x)\notag\\
 b_4(x)&=  45 (1 + x) (2 + x) \tag*{$\bigcirc$}
\end{align}
\begin{rmk}\label{remark8}
\textcolor{red}{
It is useful to note, that by construction we have
$$
(n,(k+1)!) = (n,\beta_j!)
$$
for all $k+1\in\{\beta_j,\beta_j+1,\ldots, \beta_{j-1}-1\}$.
}
\end{rmk}
\textcolor{red}{
Note that Kempner \cite{kempner2,kempner3} also introduces basic null-polynomials of the form
$$
\tilde b(x) = \frac{n}{d}\prod_{i=0}^{s(d)-1}(x-i)
$$
where $d>1$ is a divisor of $n$. If $d>1$ runs through all divisors of $n$ in decreasing order, we only list polynomials which are not multiples of polynomials that already appeared.
In the present case, when $n=90$, one obtains in this way the basic null-polynomials
$$\begin{aligned}
\tilde b_1(x) & = x(x-1)(x-2)(x-3)(x-4)(x-5) \\
\tilde b_2(x) & = 3x(x-1)(x-2)(x-3)(x-4)\\
\tilde b_3(x) & = 15x(x-1)(x-2)\\
\tilde b_4(x) & = 45x(x-1)
\end{aligned}
$$
The difference stems from the fact, that we introduced a normed null-polynomial of minimal degree by defining
$$
p(x) = \prod_{i=1}^{s(n)}(x+i),
$$
whereas Kempner uses
$$
\tilde p(x) = \prod_{i=0}^{s(n)-1}(x-i).
$$}
\eex
%
%

\textcolor{blue}{Notice that the basic null-polynomial $b_{t(n)}$ is a non-zero polynomial 
of minimal degree $q(n)$ (see, e.g.,~\cite[Theorem 8]{noebiundlorenz}). This fact is
used in the following decomposition theorem. With the notations above we have:}
\bt\label{deco}
Every null-polynomial $p$ in $\Z_n[x]$ has a unique decomposition
of the form$$p(x)=\sum_{k=1}^{t(n)}q_k(x) b_k(x),$$
where $q_k\in \Z_{n/\alpha_k}[x]$ has degree strictly
less than $\beta_{k-1}-\beta_k$ if $k>1$ and where
$\operatorname{deg}(q_1)=\operatorname{deg}(p)-\beta_1$.
\et
\proof
We start by proving the existence of a decomposition of the desired 
type.

In a first step, we can write
$$p(x)=q_1(x)b_1(x)+p_1(x)$$
with $q_1\in\Z_n[x]$, $\operatorname{deg}(q_1)=\operatorname{deg}(p)-\beta_1$,
and $\operatorname{deg}(p_1)<\beta_1$,
by dividing the polynomials with remainder (observe that $b_1$ is normed).

Now, we assume by induction that the decomposition has the form
$$
p(x)=\sum_{k=1}^lq_k(x)b_k(x)+p_l(x)
$$
with $\operatorname{deg}(p_l)<\beta_l$. Then, the next step is 
carried out as follows:
$p_l$ is a null-polynomial in $\Z_n[x]$ of the form
$$
p_l(x)=a_0+a_1x+\ldots+a_{\beta_l-1}x^{\beta_l-1}.
$$
Hence, by Theorem~\ref{voll}, it follows that
$$
a_i(\beta_l-1) ! \equiv 0\mod n
$$
for all $i\in\{0,1,\ldots,\beta_l-1\}$.
Since $\beta_{l+1}=s((n,(\beta_l-1) ! ))<\beta_l$, this implies
$$
\alpha_{l+1}\mid a_i
$$
for all $i\in\{0,1,\ldots,\beta_l-1\}$. Hence, we can divide the polynomial
$p_l$ by $b_{l+1}$ with remainder and obtain
$$p_l(x)=q_{l+1}(x)b_{l+1}(x)+p_{l+1}(x)$$
with $\operatorname{deg}(p_{l+1})<\beta_{l+1}$, 
$\operatorname{deg}(q_{l+1})<\beta_l-\beta_{l+1}$ and 
$q_{l+1}\in\Z_{n/\alpha_{l+1}}[x]$.
This iterative process ends as soon as $\operatorname{deg}(p_{l+1})<q(n)$,
since then, it follows that $p_{l+1}\equiv 0\mod n$ \textcolor{red}{by~\cite[Theorem 8]{noebiundlorenz}}.

Now, we assume by contradiction that there exist two different 
decompositions of $p$, say
\begin{equation}\label{eq-xyz}
\textcolor{red}{0\equiv \sum_{k=1}^{t(n)} b_k(q_k-\tilde q_k)\mod n}
\end{equation}
with a smallest index $k_0$ with $q_{k_0}\neq \tilde q_{k_0}$.
\textcolor{red}{Let $i$ denote  the highest power $i$ in $q_{k_0}$ and $\tilde q_{k_0}$ with
different coefficients $a_i\neq \tilde a_i$ in $\Z_{n/\alpha_{k_0}}$. Then, according to the construction
of the basic null-polynomials $b_k$, the coefficient of the highest power of $x$
on the right hand side of~(\ref{eq-xyz}) is  $\alpha_{k_0}(a_i-\tilde a_i)$. 
By~(\ref{eq-xyz}), we have
$$
\alpha_{k_0}\underbrace{(a_i-\tilde a_i)}_{\in\Z_{n/\alpha_{k_0}}}
\equiv 0\mod n
$$
which implies that $a_i\equiv \tilde a_{i}\bmod (n/\alpha_{k_0})$, and  
this is a contradiction.}
\endproof

\subsection{The number of Polyfunctions}\label{nop}
The result of the previous section allows now to compute the 
cardinality of the ring $G(\Z_n)$.
\bc\label{co1}
The number $\Psi(n)$ of polyfunctions over $\Z_n$ is given by
$$
\Psi(n)=\prod_{k=1}^{t(n)}(n,\beta_k ! )^{\beta_k-\beta_{k-1}}
$$
with the convention $\beta_0:=0$.
\ec

\proof
We consider the additive group $F(n)$ of polynomials in $\Z_n[x]$ 
of degree strictly less than $s(n)$
and the normal subgroup $N(n)$ of all null-polynomials in $F(n)$.
The additive group of polyfunctions over $\Z_n$ is
then isomorphic to the quotient $F(n)/N(n)$.
All cosets have the cardinality 
of the set of null-polynomials of degree strictly less than $s(n)$, namely,
according to Theorem~\ref{deco},
$$
|N(n)|=\prod_{i=2}^{t(n)}\left(\frac n{\alpha_i}\right)^{\beta_{i-1}-\beta_i}.
$$
On the other hand, the number of polynomials of degree strictly less than
$s(n)$ is $|F(n)|=n^{\beta_1}$. Division $|F(n)|/|N(n)|$ gives the claimed formula.
\endproof


\bex\label{counting-zero-polynomials-2} Let us come back to Example \ref{counting-zero-polynomials-1} with $n=90$:
The formula in Corollary~\ref{co1}
gives $\Psi(90)=(90,6 ! )^6(90,5 ! )^{-1}(90,3 ! )^{-2}(90, 2 ! )^{-1}=246037500$
for the number of polyfunctions over $\Z_{90}$.\hfill$\bigcirc$
\eex
In the case when $n$ equals the power of a prime number
the formula for $\Psi$ takes a particularly simple
form. Since $\Psi$ will be
shown to be multiplicative, it is actually enough to 
know the values of $\Psi(p^m)$ for $p$ prime (see the following Section~\ref{phm}).

\subsubsection{The Case \boldmath{$n=p^m$}, \boldmath{$p$} prime}\label{phm}
At this point it is useful to include a general remark on rings of
polyfunctions: If $R$ and $S$ are commutative rings with 
unit element, then $G(R\oplus S)$ 
and $G(R)\oplus G(S)$ are isomorphic as rings in the obvious way.
In particular, since $\Z_n\oplus\Z_m\cong\Z_{nm}$ if $m$ and $n$ are
relatively prime, we
have that
$$G(\Z_{nm})\cong G(\Z_n)\oplus G(\Z_m)$$ if $(m,n)=1$.
Therefore, we may confine ourselves to the case $n=p^m$, $p$ prime,
without loss of generality.

This observation gives rise to the following version of Corollary~\ref{co1}, \textcolor{red}{see also~\cite{keller}}.
\bp\label{ppp}
Let $\Psi(n)$ denote the number of polyfunctions over $\Z_n$
and $s$ the Smarandache function. Then,
\begin{itemize}
\item[(i)] the function $\Psi$ is multiplicative, i.e.\ if
$(m,n)=1$ then $\Psi(mn)=\Psi(m)\Psi(n)$, and
\item[(ii)] for a prime number  $p$ and $m\in\N$ there holds
$$
\Psi(p^m)=\exp_p\left(\sum_{k=1}^ms(p^k)\right),
$$ 
where we write $\exp_pa:=p^a$ for typographical reasons.
\end{itemize}
\ep
\bex\label{counting-zero-polynomials-3} Before we prove Proposition~\ref{ppp}, we come back
to Example \ref{counting-zero-polynomials-2}, where $n=90$.
By (i) in Proposition~\ref{ppp}, we have
$$\Psi(90)=\Psi(2)\Psi(3^2)\Psi(5)$$
and the factors are by (ii)
$\Psi(2)=2^2$, $\Psi(3^2)=3^{3+6}$ and $\Psi(5)=5^5$. The product
of these numbers is $\Psi(90)=
4\cdot 19683\cdot 3125=246037500$ in accordance with
the calculation in Example~\ref{counting-zero-polynomials-2}.\hfill$\bigcirc$
\eex
At this point it is useful to
introduce one more quantity which will play a role
in the proof of Proposition~\ref{ppp} and which is going to be used
in the description of the algebraic structure of
the ring of polyfunctions over $\Z_n$ (see Section~\ref{riri}).
For prime numbers $p$ and integers $k\ge 0$
we define
$$
e_p(k):=\max\{x\in \N_0\,:\,p^x\mid k ! \}.
$$ 
Notice that $e_p(k)=j$ for $jp\le k<(j+1)p$ if $k<p^2$. But the next number is
$e_p(p^2)=p+1$.

{\bf Proof of Proposition~\ref{ppp}}\\
(i) The multiplicativity 
follows immediately from the remark preceding the proposition.

(ii) The basic null-polynomials of degree strictly less
than $s(p^m)$ are in this case \textcolor{red}{(see~\eqref{basicnullpolynomial})} given by
$$b_k(x)=p^{m-e_p(k)}\prod_{i=1}^{k}(x-i)$$
for $k=p,2p,3p,\ldots,s(p^m)-p$.
Thus the number of null-polynomials in $\Z_{p^m}[x]$
of degree strictly less than $s(p^m)$ is 
$$
\prod_{k=1}^{\frac{s(p^m)}p-1}p^{p e_p(p k)},
$$
and the total number of polynomials in $\Z_{p^m}[x]$
of degree strictly less than $s(p^m)$ is 
$$p^{m s(p^m)}. $$
Division of both numbers yields the number of polyfunctions over 
$\Z_{p^m}$, namely
$$\Psi(p^m)=\exp_p\left(p\sum_{k=0}^{\frac{s(p^m)}p-1}(m-e_p(pk))\right).$$
Hence, the claim is proved if we verify that for all $m\in \N$ there holds
\begin{equation}\label{indu}
p\sum_{k=0}^{\frac{s(p^m)}p-1}(m-e_p(pk))=\sum_{k=1}^m s(p^{\color{red}{k}}).
\end{equation}
Obviously, (\ref{indu}) is true for $m=1$. Moreover 
$s(p^{m+1})-s(p^m)$ is either $0$ or $p$. Using this, it is easy to see,
that~(\ref{indu}) holds for $m+1$ if it is correct for $m$, and
the claim follows by induction.
\endproof
\begin{rmk}\phantom{a}
\begin{enumerate}
\item[(i)]The formula in (ii) above is particularly simple
in the case $m\le p$: 
We observe that $s(p^k)=kp$ for $k\le p$. Thus
$$\sum_{k=1}^m s(p^k)=p\binom{m+1}2\text{ and }\Psi(p^m)=\exp_p\left(p\binom{m+1}2\right)$$
for  $m\le p$.
\item[(ii)]\textcolor{black}{While the present approach for counting the number of polyfunctions in $\Z_n$ consists in finding a unique representative for each null-polynomial, in~\cite[Theorem 5, p.\ 8]{noebi}, each polyfunction is shown to have a unique representative. An alternative proof of Theorem~\ref{ppp} above is then given in~\cite[Theorem 6, p.\ 9]{noebi} by counting these representatives. Moreover, a very short formula for $\Psi(n)$ is given in~\cite[Theorem 9, p.\ 10]{noebi} in terms of the Smarandache function, the Mangoldt function and the Dirichlet convolution. 
}
\item[(iii)]\label{additive-structure}
Not only the formula for $\Psi(n)$ looks particularly
pleasant for $n=p^m$, also the decomposition of the
additive group $F(n)$ takes its simplest form for
powers of prime numbers. As mentioned
earlier in this section, it is sufficient to
know the structure of $F(n)$ for $n=p^m$. In this case,
the decomposition in Theorem~\ref{thth} below simplifies to
$$
F(p^m)\cong p\bigoplus_{k=0}^{s(p^m)/p-1} \Z_{p^{m-e_p(pk)}}.
$$
\textcolor{red}{Here and throughout Section \ref{algebra}, we will use the notation
$$
nG = \bigoplus_{i=1}^n G
$$ for the $n$-fold direct product of a group $G$ with itself, where $n\in\N$.}
\end{enumerate}
\end{rmk}
\section{Algebraic properties of the ring of polyfunctions}\label{algebra}
\subsection{The additive group of Polyfunctions}\label{addi}
Let $F(n)$ denote the additive group of polyfunctions over $\Z_n$
and $F_k(n)$ the subgroup of polyfunctions which have a representative
of degree less than or equal to $k$.
Using the notation of \textcolor{red}{Section \ref{dec}}, we have the following
result:
\bt\label{thth}
The group $F(n)$ is isomorphic to
$$
\bigoplus_{j=1}^{t(n)}(\beta_j-\beta_{j+1})\Z_{\alpha_{j+1}}
$$
with the convention $\beta_{t(n)+1}:=0$ and $\alpha_{t(n)+1}:=n$.
\et
We prepare the proof by the following lemma:
\bl\label{lele}
Let $\beta_j\le k+1<\beta_{j-1}$, $k\ge 0$, $2\le j\le t(n)+1$.
Then there holds:
\begin{itemize}
\item[(i)] Every element in the quotient
$F(n)/F_k(n)$ has order less than or equal to $\alpha_j$.
\item[(ii)] The polyfunction represented by $x^{k+1}$ has the
order $\alpha_j$ in $F(n)/F_k(n)$.
\end{itemize}
\el
{\bf Proof of the Lemma}\\
(i) We have, that in $F(n)/F_k(n)$
$$
\alpha_j x^{k+1}=\alpha_j x^{\beta_j} x^{k+1-\beta_j}=
\underbrace{b_j(x)}_{\parbox{0mm}{\parbox{10cm}{$=0$ for all $x\in
\Z_n$}}}x^{k+1-\beta_j}=0
$$
since $\beta_j\le k+1$. \textcolor{red}{Here, $b_j$ is a basic null-polynomial (see Section \ref{dec}).} Now, every $f\in F(n)/F_k(n)$ contains
$x^{k+1}$ as a factor and hence $\operatorname{ord}(f)\le\alpha_j$.

(ii) Suppose $\alpha x^{k+1}=0$ in $F(n)/F_k(n)$ for some $\alpha$ in $\Z_n$.
Then, by Theorem~\ref{voll},
$\alpha(k+1) ! \equiv 0\mod n$. Hence, $\alpha$ is a multiple of
$$
\frac{n}{(n,(k+1) ! )}>\frac{n}{(n,\beta_{j-1} ! )}=\alpha_{j-1}
$$
since $k+1<\beta_{j-1}$. Thus we have
$$\frac{n}{(n,(k+1) ! )}\ge \alpha_j$$
\textcolor{red}{(see Remark~\ref{remark8})} and hence $\alpha\notin\{1,2,\ldots,\alpha_{j}-1\}$.
\endproof
Now, Theorem~\ref{thth} follows from Lemma~\ref{lele} by iteration:
First, we observe that $1\in F(n)$ has the (maximal) order $n=\alpha_{t(n)+1}$.
Thus 
$$F(n)\cong \Z_n\oplus F(n)/F_0(n)$$
since  finite Abelian groups split off a maximal cyclic subgroup.
Now, we proceed iteratively and split in each step
$$
F(n)/F_k(n)\cong \Z_{\alpha_j}\oplus F(n)/F_{k+1}(n)
$$
by using Lemma~\ref{lele}. The process stops as soon as
$k+1=s(n)$, and by collecting the quotients we obtain the claimed decomposition.
\endproof
\bex We revisit Example~\ref{counting-zero-polynomials-1}, \ref{counting-zero-polynomials-2} and \ref{counting-zero-polynomials-3} respectively in order to compute the decomposition of $F(90)$. With the notational conventions of Theorem~\ref{thth} we have:
\begin{center}
\begin{tabular}{r|ccccc}
$j$ & 1 & 2 & 3 & 4 & 5\\
\hline
$\alpha_j$ & 1 & 3 & 15 & 45 & 90\\
$\beta_j$ &  6 & 5 & 3 & 2 & 0
\end{tabular}
\end{center}
In a first step we decompose $$F(90) \cong  \mathbb{Z}_{90}\oplus F(90)/F_0(90).$$
If $k=0$, we have $\beta_5 < k+1 < \beta_4$ and hence
$F(90)/F_0(90)$ splits off a cylic subgroup of order $\alpha_5=90$ and hence $F(90)/F_0(90) \cong \mathbb Z_{90}\oplus F(90)/F_1(90)$.\\
If $k=1$, we have $\beta_4\leqslant k+1 <\beta_3$ and hence $F(90)/F_1(90)$ splits off a cyclic subgroup of order $\alpha_4$ and hence
$F(90)/F_1(90) \cong \mathbb Z_{45}\oplus F(90)/F_2(90).$

If $k=2,3$, we have $\beta_3\leqslant k+1<\beta_2$ so we might split off twice the subgoup $\mathbb Z_{15}$ and hence
$
F(90)/F_2(90) \cong \mathbb Z_{15}\oplus\mathbb Z_{15} \oplus F(90)/F_4(90).
$

Finally, if $k=4$, it holds that $\beta_2\leqslant k+1<\beta_1$ and we find
$
F(90)/F_4(90) \cong \mathbb Z_3
$
and the process ends. This leads to the desired decomposition
$$
F(90) \cong \mathbb Z_3 \oplus 2\mathbb Z_{15}\oplus \mathbb Z_{45}\oplus 2\mathbb Z_{90}
$$
and we find again $|F(90)|=3\cdot 15^2\cdot 45\cdot 90^2 = 246037500$ in accordance with Example~\ref{counting-zero-polynomials-2} and \ref{counting-zero-polynomials-3}.
\eex
\begin{rmk}
Since it turns out that it is sufficient to know the structure of $F(p^m)$ for prime numbers $p$ (see Section~\ref{phm}), observe that in this case, the decomposition described in 
Theorem~\ref{thth} takes a particularly simple form (see Remark \ref{additive-structure}, item (iii), above).\end{rmk}

\subsection{The ring of polyfunctions}\label{riri}
In this section, we use the shorthand notation
$G(n)$ for $G(\Z_n)$, i.e.\ the ring of polyfunctions over $\Z_n$.
We recall that $G(mn)\cong G(m)\oplus G(n)$ if $(m,n)=1$, 
and hence we may restrict ourselves to
investigate the
structure of $G(n)$ in the case $n=p^m$ for $p$ prime.
Let $I_{p,m}$ be the ideal of polynomials in $\Z_{p^m}[x]$
defined by
$$
I_{p,m}=\{f\in \Z_{p^m}[x]\,:\, \text{$f(kp)=0$ for all $k$}\}.
$$
Then, we have the following decomposition:
\bt\label{ring}
\begin{itemize}
\item[(i)] $G(p^m)\cong p\,\Z_{p^m}[x]/I_{p,m}$\,.
\item[(ii)] $\Z_{p^m}[x]/I_{p,m}$ is not decomposable.
\end{itemize}
\et
\proof
We proceed in several steps:

1.\ step: For $j\in\{0,1,\ldots, p-1\}$ let
$$
R_j(p^m):=\{f\in G(p^m)\,:\, \text{$f(k)=0$ if $k\not\equiv j\mod p$}\}.
$$
It is clear that $R_j(p^m)$ is an ideal of $G(p^m)$ and that
$R_i(p^m)\cap R_j(p^m)=\{0\}$ if $i\neq j$.

2.\ step: We show that 
$G(p^m)\cong\displaystyle{\bigoplus_{j=0}^{p-1}R_j(p^m)}$.

To see this, we define
$$\epsilon_0(x):=1-x^{m\phi(p^m)},$$
where $\phi$ denotes Euler's $\phi$-function. Then we have
$$
\epsilon_0(k)\equiv\begin{cases}0&\text{if $k\not\equiv 0\mod p$}\\
1 &\text{if $k\equiv 0\mod p$}\end{cases}
\mod p^m.
$$
Moreover, for $\epsilon_j(x):=\epsilon_0(x-j)$, we have similarly
$$
\epsilon_j(k)\equiv\begin{cases}0&\text{if $k\not\equiv j\mod p$}\\
1 &\text{if $k\equiv j\mod p$}\end{cases}
\mod p^m.
$$
Hence, for $f\in G(p^m)$, we have $f \epsilon_j\in R_j(p^m)$ and
$$f=\sum_{j=0}^{p-1}f \epsilon_j.$$ Then,
$$
\Phi_0:G(p^m)\to\bigoplus_{j=0}^{p-1}R_j(p^m), f\mapsto (
f\epsilon_0,f\epsilon_1,\ldots,f\epsilon_{p-1})
$$
is a ring isomorphism (the ring operations $+$ and $\cdot$ are,
as usual,
defined component\-wise).

3.\ step: We show that $R_j(p^m)\cong R_0(p^m)$ for $j\in\{0,1,\ldots, p-1\}$.

The map
$$
\Phi_1:R_0(p^m)\to R_j(p^m), f\mapsto \color{red}{g,}
$$\textcolor{red}{where $g(x):=f(x-j)$, $x\in \Z_{p^m}$} is a ring isomorphism. Hence, according to the second step, we have that
$$
G(p^m)\cong p R_0(p^m).
$$
4.\ step: We show that $R_0(p^m)\cong \Z_{p^m}[x]/I_{p,m}$.

To see this, we consider the map
$$
\Phi_2:\Z_{p^m}[x]\to R_0(p^m), f\mapsto f \epsilon_0.
$$
$\Phi_2$ is a surjective ring homomorphism.
If $f\in\operatorname{ker}(\Phi_2)$, then $\Phi_2(f)(k)=0$ for all $k
\in\Z_{p^m}$
and hence $f(jp)\epsilon_0(jp)=f(jp)=0$ for all $j$. This implies that
$f\in I_{p,m}$. Arguing in the opposite
direction, we conclude that $f\in I_{p,m}$
implies that $f\in\operatorname{ker}(\Phi_2)$.

Now, (i) follows from the third and the 
fourth step and it remains to prove (ii).
This is done in the last step:

5.\ step: We show, that $R_0(p^m)$ is not decomposable:

Let $f\in R_0(p^m)$ be such that $f^2=f$. In particular, this means
$f^2(jp)=f(jp)$ for all $j$. Hence, $f(jp)\in\{0,1\}$ for all
$j$. Observe, that
$$
f(jp)\equiv f(0)\mod p
$$
and hence
$$ f(k)=0\qquad\text{for all $k\in\Z_{p^m}$}$$
or
$$
f(k)=\begin{cases}
0&\text{if $k\not\equiv 0\mod p$,}\\
\text{$1$}&\text{if $k\equiv 0\mod p$.}
\end{cases}
$$
It follows that only two elements $f\in R_0(p^m)$ with the property 
$f^2=f$ exist.
In a decomposable ring there are at least four elements
with $f^2=f$. This completes the proof.
\endproof

We now want to investigate the structure of the
ideal $I_{p,m}$ in more detail. 
First, for $m\in \N$ and a prime number $p$, we define
$$
s^*(p^m):=\min\{x\in\N\,:\,p^m\mid p^xx ! \}.
$$
Then, for $r\in\{1,2,\ldots,s^*(p^m)-1\}$ let
$$
e^*(r):=\max\{x\in \N\,:\, p^x\mid \color{red}{p^r r !} \}
$$
and 
$$
e^*(s^*(p^m)):=m.
$$
\begin{rmk} $s^*$ is connected with the Smarandache function by
$$
p\, s^*(p^m)=s(p^m).
$$
\end{rmk}
Let us assume, that $f\in I_{p,m}$:
$$
f(x)=a_1x+a_2x^2+\ldots+a_rx^r.
$$
Then, $f(jp)\equiv 0\mod p^m$ for all $j$ and hence, the polynomial
$$
g(x):=a_1px+a_2p^2x^2+\ldots+a_rp^rx^r
$$
is a null-polynomial over $\Z_{p^m}$. Hence, it follows from
Theorem~\ref{voll} that
$$
a_kp^kr ! \equiv 0\mod p^m
$$
for all $k\in\{1,2,\ldots,r\}$. From, this congruence, we
immediately obtain the following conclusion.
\bp\phantom{1}
\begin{itemize}
\item[(i)]
If $f\in I_{p,m}$ is normed, then $\operatorname{deg}(f)\ge s^*(p^m)$.
\item[(ii)]
If $f\in I_{p,m}$, $f(x)=a_1x+a_2x^2+\ldots+a_rx^r$, with $r\le s^*(p^m)$,
then
$$
p^{m-e^*(r)+r-k}\mid a_k
$$
holds for all $k\in\{1,2,\ldots,r\}$.
\end{itemize}
\ep
Now, the polynomials in $I_{p,m}$ can be decomposed 
similarly as the 
null-polynomials (see Section~\ref{dec} \textcolor{red}{and~\eqref{basicnullpolynomial}}). The basic polynomials
are in this case
$$
b^*_k(x):=p^{m-e^*(k)}\prod_{j=1}^k (x+jp)
$$
for $k\in\{1,2,\ldots,s^*(p^m)\}$. In fact, we have:
\bl
$b_k^*\in I_{p,m}$ for all $k\in\{1,2,\ldots,s^*(p^m)\}$.
\el
\proof
We have
\begin{eqnarray}\label{lall}\color{red}
b_k^*(ip)&=&p^{m-e^*(k)}\prod_{j=1}^k (\color{red}{ip}+jp)\nonumber\\
&=&p^{m-e^*(k)}p^k\binom{\color{red}i+k}{k}k ! 
\end{eqnarray}
The right hand side of~(\ref{lall}) is congruent 0 modulo
$p^m$ for all $j$ as is easily seen by treating  separately the
cases $k<s^*(p^m)$ and $k=s^*(p^m)$.
\endproof
%
\subsection{The units in $G(\Z_{n})$}
The previous results on the algebraic structure of the ring of
polyfunctions over $\Z_n$ allow now to answer more specific
questions. As an example, we consider the multiplicative subgroup $U_n$
of units in $G(\Z_n)$ and ask for the size of $U_{3^k}$.

For this, we consider the set $Q$ of polynomials in $\Z_{3^k}[x]$
with degree strictly less than $s(3^k)=:r+1$. A polynomial
$q\in Q$, $q(x)=a_0+a_1 x+a_2 x^2+\ldots+a_r x^r$ with
$a_i\in\Z_{3^k}$, represents \textcolor{black}{according to~\cite[Proposition 3, p.\ 5]{noebi}
an invertible polyfunction (i.e.\ a unit in $G(\Z_{3^k})$) if
and only if its image is contained in the multiplicative subgroup of units in $\Z_{3^k}$, that is}
\begin{equation}\label{drei}
q(i)\not\equiv 0\mod 3\text{ 
\,\,for $i=0,1,2$.}
\end{equation} 
(Observe that $q(x+3j)\equiv q(x)\mod 3$ for all
integers $x$ and $j$.)
Let 
$$\Sigma_1:=\sum_{\substack {i=1\\ \text{$i$ odd}}}^r a_i$$
and
$$\Sigma_2:=\sum_{\substack {\color{black}{i=2}\color{black}\\\text{$i$ even}}}^r a_i.$$
Then, we can rewrite~(\ref{drei}) in the form
\begin{equation}\label{easy}\left.
\begin{array}{rcl}
a_0&\not\equiv& 0\mod 3\\
a_0+\Sigma_1+\Sigma_2 &\not\equiv& 0\mod 3\\
a_0+\Sigma_1+2\Sigma_2 &\not\equiv& 0\mod 3
\end{array}\right\}
\end{equation}
It is then easy to determine the total number $X$ of 
solutions $(a_0,a_1,\ldots,a_r)\in\Z_{3^k}^{r+1}$ of~(\ref{easy}):
$$X= 8\cdot 3^{k(r+1)-3}.$$
Now, two polynomials in $Q$  represent the same unit
in $G(\Z_{3^k})$ if and only if their difference is a null-polynomial
of degree strictly less than $s(3^k)$. The number $Y$ of such null-polynomials
is according to Proposition~\ref{ppp} given by
$$
Y=\frac{3^{ks(3^k)}}{\Psi(3^k)}.
$$
Division of $X$ by $Y$ yields  the  following result:
\begin{pro}\label{units}
$$
|U_{3^k}| = \Big(\frac 23\Big)^3\,\Psi(3^k)=\Big(\frac 23\Big)^3\,
\exp_3\left(\sum_{i=1}^ks(3^i)\right).
$$
In other words, the fraction of units among all polyfunctions
in $G(\Z_{3^k})$ is $\frac 8{27}$, independently of $k$.
\end{pro}
Proposition~(\ref{units}) gives a flavour of a more general result: 
In Section~\ref{micha} we will determine the number
of units in the ring $G_d(\mathbb Z_{p^m})$ of multivariate polyfunctions.

\section{Polyfunctions in several variables}\label{xy}
In order to keep the formulas short, we use 
the following multi-index notation: 
For $\k=
( k_1,k_2,\ldots,k_d )\in\N_0^d$ and $\x:=( x_1,x_2,\ldots,x_d)\color{red}{\in\N_0^d}$
let
$$\x^\k:=\prod_{i=1}^dx_i^{k_i},\quad
\k ! :=\prod_{i=1}^dk_i !,\quad
|\k|:=\sum_{i=1}^d k_i,\color{red}{\text{\quad and\quad}\binom{\x}{\k}}:=\prod_{i=1}^d \,\binom{x_i}{k_i}.
$$%
Recall that
$$G_d(R)=\{f:R^d\to R\mid 
\exists p\in R[x_1,x_2,\ldots,x_d]\,\, \forall x\in R^j\implies p(x)=f(x)\},$$ 
equipped with pointwise addition and multiplication denotes the ring of polyfunctions in $d$ variables, whenever $R$ is a commutative ring with unit element.

An alternative (but equivalent) construction is to
define $G_d(R)$ recursively as the ring of polyfunctions in
one variable
from $R$ to $G_{d-1}(R)$ by
$$ G_d(R)=\{f:R\to G_{d-1}(R)\mid 
\exists p\in G_{d-1}(R)[x]\,\, \forall x\in R\implies p(x)=f(x)\}.
$$
\subsection{The number of multivariate polyfunctions on $\mathbb Z_n$}

%
%
We recall a few facts and definitions from~\cite{noebi} in order to count the number of polyfunctions on $\Z_n$ in $d$ variables, and again it is enough to find a formula for $n=p^m$ since we have the natural decomposition $G_d(\Z_{ab})\cong G_d(\Z_{a})\oplus G_d(\Z_{b})$ if $(a,b)=1$. We define the set
\begin{equation}\label{sd(n)}
\color{red}S_d(n) := \{\k\in\N_0^d: n \nmid \k!\}
\end{equation}
and let $s_d(n) :=|S_d(n)|$ be the generalization of the Smarandache function introduced in~\cite{noebi}. As for the case of one variable we define
$$
e_p(\k):=\max\{x\in\N_0\,:\, p^x\mid \k!  \}.
$$
\begin{defi}
Let $a$ be an element of $\Z_n$.
We say, the polynomial $a\x^\k\in \Z_n[\x]$ is {\em reducible\/} (modulo $n$)
if a polynomial $p(\x)\in \Z_n[\x]$ exists with $\operatorname{deg}(p)
< |\k|$ such that $a\x^\k\equiv p(\x)\mod n$ for all $\x\in \Z_n^d$.
Moreover, we say that $a\x^\k$ is {\em weakly reducible\/}
if $a\x^\k\equiv p(\x)\mod n$ for all $\x\in \Z_n^d$, where $p\in \Z_n[\x]$ 
is such that $\operatorname{deg}(p)\le |\k|$ (instead of $\operatorname{deg}(p)< |\k|$)
and such that $\x^\k$ (or a multiple of it) does not appear as a monomial in $p$.
\end{defi}
 We will need the following lemma (see also~\cite[Lemma 4, p.\ 6]{noebi}) which characterizes tuples $\k$ for which $a\x^\k$ is (weakly) reducible in $\Z_n[\x]$.
\bl\label{lll}\phantom{1} \begin{itemize}
\item[(i)] If $a\x^\k$ is weakly reducible modulo $n$, then $n\mid a\k!$.
\item[(ii)] If $n\mid a\k!$, then $a\x^\k$ is reducible modulo $n$.
\end{itemize}
\el
\begin{proof}
(i)
We assume, that $p(\x)$
reduces $a\x^\k$ weakly. Hence,
$q(\x):=a\x^\k-p(\x)$ is a null-polynomial in $d$ variables over $\Z_n$.
Let us define the following ``integral'' for functions $f:\Z_n\to \Z_n$:
$$
\int_0^m f(x)d\mu(x):=\sum_{j=0}^m(-1)^{m-j}\binom mj f(j).
$$
Now, we write $q$ in the form
$$
q(\x)=\sum_{\substack{\l\in\N_0^d\\|\l|\le|\k|}}q_{\l}\x^\l
$$
for suitable coefficients $q_\l\in\Z_n$, with $q_\k=a$.
Then, modulo $n$, we have
\begin{eqnarray*}
0&=&\int_0^{k_d}\int_0^{k_{d-1}}\dots\int_0^{k_1}q(\x)d\mu(x_1)\ldots d\mu(x_{d-1})d\mu(x_d)=\\
&=&\sum_{\substack{\l\in\N_0^d\\|\l|\le|\k|}}q_\l\int_0^{k_d}\int_0^{k_{d-1}}\dots\int_0^{k_1}
\x^\l d\mu(x_1)\ldots d\mu(x_{d-1})d\mu(x_d).
\end{eqnarray*}
Observe that the only term which does not vanish in the above sum is
$$
q_\k\int_0^{k_d}\int_0^{k_{d-1}}\dots\int_0^{k_1}
\x^\k d\mu(x_1)\ldots d\mu(x_{d-1})d\mu(x_d)=a\k ! .
$$
In fact all other terms vanish by~(\ref{special}), since $|\l|\le|\k|$ and
$\l\neq\k$ implies that for some $i\in\{0,1,\ldots,d\}$ we have $l_i< k_i$ and
therefore the integral with respect to $x_i$ gives zero.
This completes the proof of (i).

(ii) We assume, that $n\mid a\k !  $.
Then, the polynomial
$$
q(\x):=a\prod_{i=1}^d\,\prod_{l=1}^{k_i}(x_i+l)=a \k! \prod_{i=1}^d \,\binom{x_i+k_i}{k_i} = a\k !\binom{\x + \k}{\k}
$$
is a  null-polynomial over $\Z_n$ and the term of maximal degree
is $a\x^\k$. Hence, $q(\x)-a\x^\k$ reduces $a\x^\k$.
\end{proof}
As an immediate consequence, we have:
\begin{cor}
A monomial $\x^\k$ is reducible modulo $n$ if and only if it is weakly reducible. 
\end{cor}
Furthermore it is proved in~\cite[Proposition 5, p.\ 8]{noebi} that every polyfunction $f\in G_d(\Z_{p^m})$ has a unique representative of the form
$$
f(\x) = \sum_{\k\in\N_0^d \atop e_p(\k)<m}\alpha_\k \x^\k,
$$
where $\alpha_\k\in\{0,1,\ldots,p^{m-e_p(\k)}-1\}$. Notice, that $e_p(\k)<m$ if and only if $\k\in S_d(p^m)$ and hence this representative can be written as
\begin{equation}\label{representative}
f(\x) = \sum_{\k\in S_d(p^m)}\alpha_\k \x^\k.
\end{equation} In the case of one variable, what the Smarandache function
really does is counting the number of monomials $x^k$, $k\in \N_0$,
which are not reducible. Using the unique representative of a polyfunction above we can count
the number of monomials $\x^\k$, $\k\in \N_0^d$, which are 
not reducible and hence to find a formula for $\Psi_d(n)$ which counts the number of polyfunctions in $G_d(\Z_n)$. \textcolor{red}{In view of \eqref{representative} we have for every coefficient $\alpha_\k$ exactly $p^{m-e_p(\k)}$ choices and therefore we obtain:}
\begin{pro}\label{alternative}
The number of polyfunctions in $G_d(\Z_{\color{red}{p^m}}\color{black})$ is given by
\begin{equation}\label{neu}
\Psi_d(p^m)=\prod_{\substack{\k\in\N_0^d\\ e_p(\k)< m}} p^{m-e_p(\k)}.
\end{equation}
\end{pro}
On the other hand it is shown in~\cite[Theorem 6, p.\ 9]{noebi} that
\begin{equation}\label{alt}
\Psi_d(p^m) = \exp_p\left(\sum_{k=1}^m s_d(p^k)\right).
\end{equation}
The equivalence of the two formulas~(\ref{neu}) and~(\ref{alt}) can be established by a similar induction argument as in the proof of Proposition~\ref{ppp}.
However, it is much more instructive, to give a direct algebraic argument:
We consider the surjective homomorphism $H$ of rings defined by
\begin{equation}\label{H}
H:G_d(\Z_{p^{m+1}})\to G_d(\Z_{p^{m}}),\quad f\mapsto H(f):=h\circ f\circ h^*.
\end{equation}
Here, 
$$ h:\Z_{p^{m+1}} \to \Z_{p^{m}},\quad [x]_{p^{m+1}}\mapsto [x]_{p^m},$$
where $[x]_n$ denotes the coset of $x\in\Z$ modulo $n$.  \textcolor{red}{Similarly,
$$ h^*:\Z^d_{p^{m}} \to \Z^d_{p^{m+1}},\quad [\x]_{p^{m}}\mapsto [\x]_{p^{m+1}}$$
where $[\x]_{p^m}=[(x_1,\ldots x_d)]_{p^m}:=([x_1]_{p^m},\ldots,[x_d]_{p^m})$
for $0\le x_i<p^m$.}
Then, 
$$
\Psi_d(p^{m+1})=|G_d(\Z_{p^{m+1}})|=|G_d(\Z_{p^m})||\operatorname{ker} H|
$$
and the equivalence of~(\ref{neu}) and~(\ref{alt}) is proved if we can show that
\begin{equation}\label{last}
|\operatorname{ker}H|=p^{s_d(p^{m+1})}.
\end{equation}

\textcolor{red}{In view of~\eqref{representative}, every polyfunction $f\in G_d(\Z_{p^{m+1}})$ has a unique representation
$$
f(\x) = \sum_{\k\in S_d(p^{m+1})}\alpha_{\k}\x^{\k},
$$
where $\alpha_{\k}\in\{0,1,\ldots,p^{m+1-e_p(\k)}-1\}$. Since every number in this set can be written in a unique way as
$$
\alpha_{\k}=\sum_{\{i\leqslant m+1: \k\in S_d(p^i)\}}p^{m+1-i}\alpha_{\k i},
$$
where $\alpha_{\k i}\in \Z_p$, all coefficients can be described as $i p^{m+1-e_p(\k)}$, $\k\in S_d(p^{m+1})$ and $i=0,1,\ldots,p-1$ (see also~\cite[Proposition 5, p.\ 8]{noebi}).
}

Observe, that $f\in\operatorname{ker}H$ if and only if
$f(\x)\equiv 0\mod p^m$, i.e.\ exactly if $pf$ vanishes as a function \textcolor{red}{$\Z_{p^{m+1}}^d\to\Z_{p^{m+1}}$}.

Now, for each $\k\in  \color{red}{S_d(p^{m+1})}$
and every $a_i:=i p^{m-e_p(\k)}$, $i=0,1,\ldots,p-1$,
the monomial $a_ip\x^\k$ is reducible modulo $p^{m+1}$ \textcolor{black}{by Lemma~\ref{lll} since $p^{m+1}\mid a_ip\k!$}. \textcolor{red}{This implies that $p^m\mid a_i\k!$ and hence the monomial $a_i\x^{\k}$ is reducible modulo $p^m$, i.e.}\ there exists a polynomial $q_{i,\k}(\x)$ of degree
strictly less than $|\k|$ which agrees modulo $p^m$ with
$a_i\x^\k$. Thus, $a_i\x^\k-q_{i,\k}(\x)$ 
represent polyfunctions in $\operatorname{ker}H$.
\textcolor{red}{By the considerations above, every $f\in \operatorname{ker}H$
has therefore a unique representation of the form
$$
f(\x)=\sum_{\k\in S_d(p^{m+1})}a_i\x^\k-q_{i,\k}(\x), \quad i\in\{0,1,\ldots,p-1\}
$$
and }hence $|\operatorname{ker}H|=p^{|S_d(p^{m+1})|}=p^{s_d(p^{m+1})}$, as
claimed.\endproof
\textcolor{red}{
\subsection{The number of units in $G_d(\Z_{p^m})$}\label{micha}
We end this discussion by coming back to the question of units in
the ring of polyfunctions (see Proposition~\ref{units}).
We denote by $U_{p^m}^d$ the multiplicative subgroup of units in $G_d(\Z_{p^m})$ and continue to use the notation $\Psi_d(p^m)=|G_d(\Z_{p^m})|$. 
We refer here to the formula 
$\Psi_d(p^m)=\exp_p(\sum_{k=1}^m s_d(p^k))$ from Proposition \ref{alternative}
, where $s_d(p^k)$ is defined in \eqref{sd(n)}. Then the following proposition holds
\begin{pro}\label{units-general}
$$
|U_{p^{m}}^d|=\left(\frac{p-1}{p}\right)^{dp}\Psi_{d}(p^{m}).
$$
\end{pro}
\begin{proof}
Using \cite[Proposition 3, p.5]{noebi} we know that the elements in $U_{p^m}^d$ are precisely the unit-valued polyfunctions in $G_d(\Z_{p^m})$.
Note that every function $\Z_p^d\to\Z$ is a polyfunction hence $|G_d(\Z_p)|=p^{dp}$ and since there are $p-1$ units in $\Z_p$, we have $$|U_p^d| = (p-1)^{dp}.$$
We use again the map
$$H:G_d(\Z_{p^{m+1}})\to G_d(\Z_{p^m}),\quad f\mapsto H(f) =  h\circ f \circ h^*$$
as defined after~(\ref{H}).
Now
$$
f\in U_{p^{m+1}}^d \Longleftrightarrow H(f)\in U_{p^m}^d.
$$
Indeed, $f\in U_{p^{m+1}}$ if and only if $f\circ h^*$ is unit valued with values in $\Z_{p^{m+1}}$ if and only if $((f\circ h^*)(\x),p^{m+1})=1$ if and only if $(H(f)(\x),p^{m+1})=1$ if and only if $(H(f)(\x),p^{m})=1$ (see also \cite[Remark 12.1, Lemma 7 and Lemma 8]{almaktry}). We conclude that
$$
|U_{p^{m+1}}^d| = |\operatorname{ker} H||U_{p^m}^d|
$$
and it follows from the proof of Proposition \ref{alternative} that $|\ker H|=p^{s_d(p^{m+1})}$.
So, inductively
$$
 |U_{p^{m}}^d| = \prod_{i=2}^{m}p^{s_d(p^{i})} |U_p^d|
 $$
and since $|U_p^d| = (p-1)^{dp}$ we find using the formula for $\Psi_d(p^m)$ of Proposition \ref{alternative}
$$
|U_{p^{m}}^d|=(p-1)^{dp} \exp_p\left(\sum\limits_{i=2}^{m}s_d(p^i)\right)=\left(\frac{p-1}{p}\right)^{dp}\Psi_{d}(p^{m}).
$$
\end{proof}}

\textcolor{red}{
\section*{Acknowledgement}
We would like to thank the referee for his or her very careful reading and for all the valuable remarks and suggestions which greatly helped to improve the quality and readability of this article.}

\end{document}